\documentclass[11pt]{article}
\usepackage[T1]{fontenc}
\usepackage[utf8]{inputenc}
\usepackage{amsmath,amssymb,amsthm,booktabs}
\usepackage[margin=1in]{geometry}
\usepackage{tikz}
\usetikzlibrary{arrows.meta,calc,positioning}
\usepackage{hyperref,xurl}
\hypersetup{colorlinks=true,linkcolor=blue!50!black,citecolor=blue!50!black,urlcolor=blue!50!black}
\newtheorem{theorem}{Theorem}[section]
\newtheorem{lemma}[theorem]{Lemma}
\newtheorem{proposition}[theorem]{Proposition}

\theoremstyle{remark}

\newcommand{\HH}{\mathbb H}
\newcommand{\Vol}{\operatorname{Vol}}
\newcommand{\ip}[2]{\langle #1,#2\rangle}
\newcommand{\OCT}{\mathcal O}
\newcommand{\DBL}{\mathcal D}
\newcommand{\CELL}{\mathcal C}
\title{On the second smallest ideal right-angled\linebreak hyperbolic four-polytope}
\author{Alexander Kolpakov}
\date{28 October 2026}
\begin{document}
\maketitle
\begin{abstract}
The regular ideal right-angled 24-cell has volume $4\pi^2/3$ and uniquely
minimizes volume among entirely ideal right-angled hyperbolic
four-polytopes. We prove that the second distinct volume is $8\pi^2/3$,
attained uniquely by the reflection double of the 24-cell across an
octahedral facet. We also show that an entirely ideal right-angled
four-polytope whose facets are regular ideal octahedra or their
doubles across triangular faces decomposes
into a chain of regular 24-cells, with opposite cutting facets in each
intermediate cell.
\end{abstract}

\section{Introduction}

The volume spectrum of right-angled hyperbolic four-polytopes connects
three different kinds of information: the incidence structure of the
polytope, the intrinsic geometry of its facets, and the shapes of its
cusps. Gauss--Bonnet expresses the volume in terms of face numbers,
but this expression alone does not determine which values occur or
which polytopes attain them. Even a short gap in the volume spectrum
therefore raises a geometric recognition problem.

Two basic extremal problems ask for the least facet number and the
least volume of a right-angled four-polytope. In the compact case,
Potyagailo and Vinberg~\cite[p.~66]{PV} asked whether the least facet
number is $120$. Guo and Zheng~\cite[Theorem~1.1 and Corollary~1.2]{GuoZheng}
answer this question affirmatively and establish the corresponding
volume minimum: the right-angled 120-cell uniquely attains both
$120$ facets and volume $34\pi^2/3$.

For entirely ideal polytopes, both analogous extremal problems were
settled in~\cite[Theorems~3 and~4]{Kolpakov}: the
regular ideal right-angled 24-cell uniquely attains the least facet
number $24$ and the least volume $4\pi^2/3$. The present paper
continues this ideal-volume problem by determining the next distinct
volume and its equality case. The natural candidate is the reflection
double of the 24-cell across one octahedral facet, which is again
entirely ideal and right-angled. We prove that its volume $8\pi^2/3$
is the second distinct volume and that it is the unique polytope
attaining that value. The exclusion of intervening values uses a
facet-defect inequality; the identification of the equality case
requires a separate argument constructing an octahedral
cutting section.

That cutting argument also gives a structure theorem beyond the
second-volume problem. When every facet is an ideal octahedron or
the reflection double of one across a triangular face, the whole
polytope decomposes into a chain of regular 24-cells. This additional
facet hypothesis is essential to the statement. We establish it from
the volume assumption before applying the structure theorem to the
equality case.

\subsection{Main results}

A polytope in this paper is convex, finite-sided, four-dimensional,
and of finite hyperbolic volume. It is \emph{entirely ideal} if every
vertex lies at infinity, and \emph{right-angled} if every dihedral
angle is $\pi/2$. Write $\CELL$ for the regular ideal right-angled
24-cell, $\OCT$ for its regular ideal octahedral facet, and $\DBL$
for the union of $\OCT$ with its reflection across a triangular
face. The polyhedron $\DBL$ has nine vertices, eight triangular
faces and three quadrilateral faces.

\begin{theorem}[Minimum volume]\label{thm:minimum}
For every entirely ideal right-angled hyperbolic four-polytope $P$,
\[
 \Vol(P)\ge\frac{4\pi^2}{3},
\]
with equality exactly when $P$ is isometric to $\CELL$.
\end{theorem}
See~\cite[Theorem~3]{Kolpakov}. We recall the
volume and incidence formulas underlying it below. Our later
all-octahedral recognition argument also recovers its equality
statement once those formulas are given.

\begin{theorem}[The second volume and its equality case]\label{thm:gap}
If $P$ is entirely ideal and right-angled and $P\not\cong\CELL$, then
\[
 \Vol(P)\ge\frac{8\pi^2}{3}.
\]
Equality holds exactly for the union of $\CELL$ and its reflection
across an octahedral facet. Its face vector, in increasing dimension,
is $(42,168,164,38)$.
\end{theorem}

We also consider polytopes satisfying the facet hypothesis
\begin{equation}\label{eq:spectrum}
 \text{every facet of $P$ is isometric to $\OCT$ or $\DBL$.}
\end{equation}
By ideal-polyhedron rigidity in dimension three~\cite{Rivin}, this
is equivalent to prescribing the corresponding combinatorial types.

\begin{theorem}[Octahedral cuts and chains]\label{thm:chain}
Suppose that $P$ satisfies \eqref{eq:spectrum}. Finitely many pairwise
disjoint regular ideal octahedral sections cut $P$ into regular ideal
right-angled 24-cells. If $P$ has $N$ facets of type $\DBL$, then
$N=8k$ for an integer $k\ge0$, and there are $k+1$ cells. Their
adjacency graph is a path. Two cuts incident to one cell use opposite
facets. Within this class the isometry type of an existing polytope
is determined by its number of cells, and
\[
 \Vol(P)=(k+1)\frac{4\pi^2}{3}.
\]
\end{theorem}
\subsection{Outline of the proof}

The proof has three stages.
\begin{enumerate}
\item The volume formula defines a nonnegative ridge defect $W$,
divisible by three. We associate to each facet two disjoint lists
of other positive-defect facets: neighbors across nontriangular
ridges, and opposite facets at nonsquare cusps. Their total defect
excludes the three possible volumes strictly between the minimum
and twice the minimum.
\item At twice the minimum, a further defect argument uses the
classification of small three-dimensional facets and their cusp
shapes to show that exactly eight facets are doubled octahedra;
every other facet is octahedral.
\item Under this facet hypothesis, distinguished triangles in the
doubled facets form closed surfaces. Incidence counting gives eight
triangles and six cusps in each component. Equal horospherical
widths give identities among the Lorentz normals that place the
triangles in one hyperplane. Its intersection with the polytope
is a regular ideal octahedron.
Cutting reduces the number of doubled facets by eight. Recognition
of the final pieces and an opposition argument give the chain
theorem. With eight doubled facets there is just one cut, so the
equality case consists of two reflected 24-cells.
\end{enumerate}
The last section supplies an exact halfspace model of the double.
The appendix describes the finite supporting algorithms, their
completeness, and the supplementary verification.

\subsection{Previous work}
The incidence and volume formulas used below are given
in~\cite[Lemmas~1--2]{Kolpakov}. By~\cite[Proposition~5]{Kolpakov},
an ideal right-angled three-polyhedron with a $k$-gonal face has
at least $2k$ vertices, with equality exactly for the $k$-antiprism.
The unique admissible nine-vertex case is given
in~\cite[Section~4.3, Case~D and Figure~6]{Kolpakov}. The equality
reduction also uses the ten-vertex octahedrite
in~\cite[Figure~7]{Kolpakov}: among the
ten-vertex types with only triangles and quadrilaterals it is the
unique admissible one, and the edge-adjacency graph of its four
quadrilaterals consists of two disjoint edges. The pentagonal
antiprism is treated separately. Rivin's rigidity theorem for ideal
three-polyhedra~\cite{Rivin} identifies these combinatorial models
geometrically once their dihedral angles are prescribed. The
all-wall acuteness property used below is the standard one for
nonobtuse hyperbolic polytopes~\cite[\S2.1]{Vinberg}.

\section{Discrete volume and facet defect}\label{sec:defect}
Write $v,e,r,f$ for the numbers of vertices, edges, ridges and facets of
$P$. An ideal vertex link is a Euclidean rectangular parallelepiped. Counting
its vertices and edges and using Euler's identity gives
\[
 e=4v,\qquad r=3v+f,\qquad
 \sum_{R\in\Omega_2(P)} k(R)=12v,
\]
where $k(R)$ is the number of vertices of a ridge $R$. Define
\[
 W(P)=\sum_R(k(R)-3).
\]
It follows that
\begin{equation}\label{eq:discrete}
 W(P)=3(v-f)\in3\mathbb Z_{\ge0},\qquad
 \Vol(P)=\frac{\pi^2}{3}(v-f+4)
 =\frac{4\pi^2}{3}+\frac{\pi^2}{9}W(P).
\end{equation}
The latter identity is the volume formula in~\cite{Kolpakov}; equivalently,
Gauss--Bonnet applied to the reflection orbifold uses
$\chi=1-f/2+r/4-e/8=(v-f+4)/4$.

Every facet $F$ is an ideal right-angled three-polyhedron. Its ideal
vertices have valency four, so, writing $v_F$ for their number, it has
$2v_F$ edges and $v_F+2$ polygonal faces. Its defect is
\begin{equation}\label{eq:facet}
 d(F):=v_F-6=\sum_{R\subset F}(k(R)-3),\qquad
 \sum_F d(F)=2W(P).
\end{equation}
Defect zero characterizes the regular ideal octahedron. A positive
defect facet has a nontriangular face of size $k\ge4$; the antiprism
bound $v_F\ge2k$ implies $d(F)\ge2$. In particular defect one does not
occur.

\section{Opposite cusp facets}\label{sec:injection}
For a facet $F$, let $q(F)$ be its number of nontriangular faces and
$s(F)$ the number of its cusps whose intrinsic horospherical rectangle
is nonsquare. The latter condition does not depend on horosphere scale.

\begin{lemma}[Opposite-cusp injection]\label{lem:injection}
Suppose the smallest positive facet defect of $P$ is $m$, and $F$ has
positive defect. There are $q(F)+s(F)$ pairwise distinct positive-defect
facets other than $F$: its neighbors across nontriangular ridges, and
the facets opposite $F$ at its nonsquare cusps. Therefore
\begin{equation}\label{eq:budget}
 2W(P)\ge d(F)+m\bigl(q(F)+s(F)\bigr).
\end{equation}
In particular, when $d(F)=m$,
$2W(P)\ge m(1+q(F)+s(F))$.
\end{lemma}
\begin{proof}
A neighbor across a nontriangle contains that same nontriangle and is
not octahedral. Distinct ridges of $F$ have distinct neighboring facets,
since the intersection of two facets is a face.

At a cusp $x$ of $F$, the face of the rectangular parallelepiped opposite
the $F$-face belongs to a facet $F_x$. The two rectangles are congruent.
If the rectangle of $F$ is nonsquare, so is that of $F_x$, which cannot
be the regular ideal octahedron.

The supporting hyperplanes of $F$ and $F_x$ are asymptotic and meet
at exactly the ideal point $x$. Hence the same $F_x$ cannot occur at
a second cusp of $F$, nor can it meet $F$ across any ridge. Thus all
facets in the two lists are distinct and differ from $F$. Each has
defect at least $m$. Summing their defects in \eqref{eq:facet} proves
the assertion.
\end{proof}

If $P\not\cong\CELL$, a positive-defect facet exists: otherwise $W=0$
in \eqref{eq:facet}, and the volume formula and the uniqueness
statement in~\cite[Theorem~3]{Kolpakov} identify $P$ with $\CELL$.

\section{The two cusp-shape calculations}
\begin{lemma}\label{lem:antiprism}
The cusp rectangles of the ideal right-angled $k$-antiprism have
aspect ratio $2\cos(\pi/k)$. Thus for $k\ge4$ all $2k$ cusps
are nonsquare.
\end{lemma}
\begin{proof}
Set $\theta=\pi/k$, $c=\cos\theta$, and take positive square roots in
\[
 h^2=\frac{1-c}{1+c},\quad r^2=\frac{2c}{1+c},\quad
 b^2=\frac{2}{c(1+c)},\quad \mu^2=\frac{2(1-c)}{c}.
\]
Use Lorentz form $(-1,1,1,1)$. The two vertex rings have heights
$\pm h$, radius $r$, and angular offset $\theta$. At
$u=(1,r\cos\theta,r\sin\theta,h)$, the four incident normals are
\[
 N=(h,0,0,1),\quad M_0=(1,b,0,h),\quad
 M_2=(1,b\cos2\theta,b\sin2\theta,h),\quad
 L=(1,b\cos\theta,b\sin\theta,-h).
\]
Their squares are $r^2,\mu^2,\mu^2,\mu^2$. Adjacent products vanish;
the opposite pairs are parallel. The identities $bc=r$ and
$\mu=2h/r$ give
\[
 \frac{M_0+M_2}{\mu}=\frac2\mu u,\qquad
 \frac Nr+\frac L\mu=\left(\frac hr+\frac1\mu\right)u.
\]
For two opposite unit normals with sum $\kappa u$, their separation
on a fixed horosphere is proportional to $\kappa$. To see this,
send $u$ to infinity in the upper-half-space model; the mirrors become
parallel vertical planes with unit horizontal normals. The two
separations therefore have ratio
\[
 \frac{2/\mu}{h/r+1/\mu}
 =\frac{2r^2}{1+h^2}=2c.
\]
The antiprism symmetries are vertex-transitive. Finally $2c>1$ for
$k\ge4$ (and equals one for $k=3$).
\end{proof}

\begin{lemma}\label{lem:nine}
An ideal right-angled three-polyhedron with nine vertices has
$d=3$, $q=3$ and $s=3$.
\end{lemma}
\begin{proof}
The antiprism bound excludes pentagons and larger faces. The complete
published nine-vertex octahedrite list consists of one admissible type:
two ideal regular octahedra glued across a triangular face
\cite[Section~4.3, Case~D2 and Figure~6]{Kolpakov}. Rigidity makes
this identification geometric. The three gluing cusps have rectangles
formed by joining two squares, hence aspect two. The other six cusps
remain square. Its three nontriangular faces are quadrilaterals.
\end{proof}

\section{The volume gap}\label{sec:gap}
\begin{lemma}\label{lem:two}
A positive-defect facet has at least two nontriangular faces. If it
has exactly two, it is an antiprism.
\end{lemma}
\begin{proof}
If its only nontriangle has size $k$, then $d=k-3$, whereas the
antiprism bound gives $d\ge2(k-3)$, a contradiction. If the two sizes
are $k\ge l\ge4$, then
\[
 d=(k-3)+(l-3)\ge2(k-3).
\]
Thus $k=l$ and equality holds in the antiprism bound; its equality
statement identifies the facet with $A_k$.
\end{proof}

\begin{proof}[Proof of the lower bound in Theorem~\ref{thm:gap}]
Choose a facet $F$ of smallest positive defect $m\ge2$.
Lemma~\ref{lem:injection} gives $2W\ge m(1+q+s)$.
The following cases are exhaustive.
\begin{enumerate}
\item If $F=A_k$, then $m=2k-6$, $q=2$, and
$s=2k=m+6$ by Lemma~\ref{lem:antiprism}. Hence
$2W\ge m(m+9)\ge22$.
\item A non-antiprism with $m=2$ is impossible: $v_F=8$ gives
$k\le4$ for every face, so its defect is two quadrilaterals and
Lemma~\ref{lem:two} applies.
\item If $m=3$, Lemma~\ref{lem:nine} gives $2W\ge3(1+3+3)=21$.
\item If $m=4$, every face has size at most five. A pentagon gives
equality in the antiprism bound and has already been treated. Otherwise
there are four quadrilaterals, so $2W\ge4(1+4)=20$.
\item If $m\ge5$ and $F$ is not an antiprism,
Lemma~\ref{lem:two} gives $q\ge3$, whence $2W\ge4m\ge20$.
\end{enumerate}
Thus $W\ge10$. Since $W$ is divisible by three, $W\ge12$.
Equation~\eqref{eq:discrete} proves the lower bound. In particular the
three prospective intermediate defects $3,6,9$, with respective
volumes $5\pi^2/3$, $2\pi^2$ and $7\pi^2/3$, cannot occur.
\end{proof}

\section{The intrinsic doubled octahedron}\label{sec:doubled}

The three vertices of the gluing face of $\DBL$ have a simple
intrinsic description: their horospherical rectangles are nonsquare.
At each of these vertices the rectangle is obtained by adjoining
two unit squares along one side, and therefore has aspect ratio two.
The remaining six rectangles are square. The three quadrilateral
faces of $\DBL$ meet pairwise at its three nonsquare cusps and have
no common edges. At a nonsquare cusp, the two quadrilateral sides
are the two long sides of its $2$ by $1$ rectangle.

\begin{figure}[htbp]
\centering
\begin{tikzpicture}[scale=1.25,>=Latex]
 \fill[blue!7] (0,0) rectangle (4,2);
 \draw[very thick,blue!65!black] (0,0)--(4,0) (0,2)--(4,2);
 \draw[thick] (0,0)--(0,2) (4,0)--(4,2);
 \draw[dashed,very thick,red!70!black] (2,0)--(2,2);
 \draw[<->] (0,-.45)--node[below] {$2$}(4,-.45);
 \draw[<->] (4.5,0)--node[right] {$1$}(4.5,2);
 \node at (1,1) {one square};
 \node at (3,1) {one square};
 \node[above,blue!65!black] at (2,2.15) {quadrilateral sides};
 \node[align=left,anchor=west,red!70!black] at (5.2,1)
 {$T_F$ meets the link\\along the middle segment};
\end{tikzpicture}
\caption{The cusp rectangle of a doubled octahedron at a gluing
vertex. The blue sides belong to its quadrilateral faces. The red
segment is the link of the internal triangle.}
\label{fig:rectangle}
\end{figure}
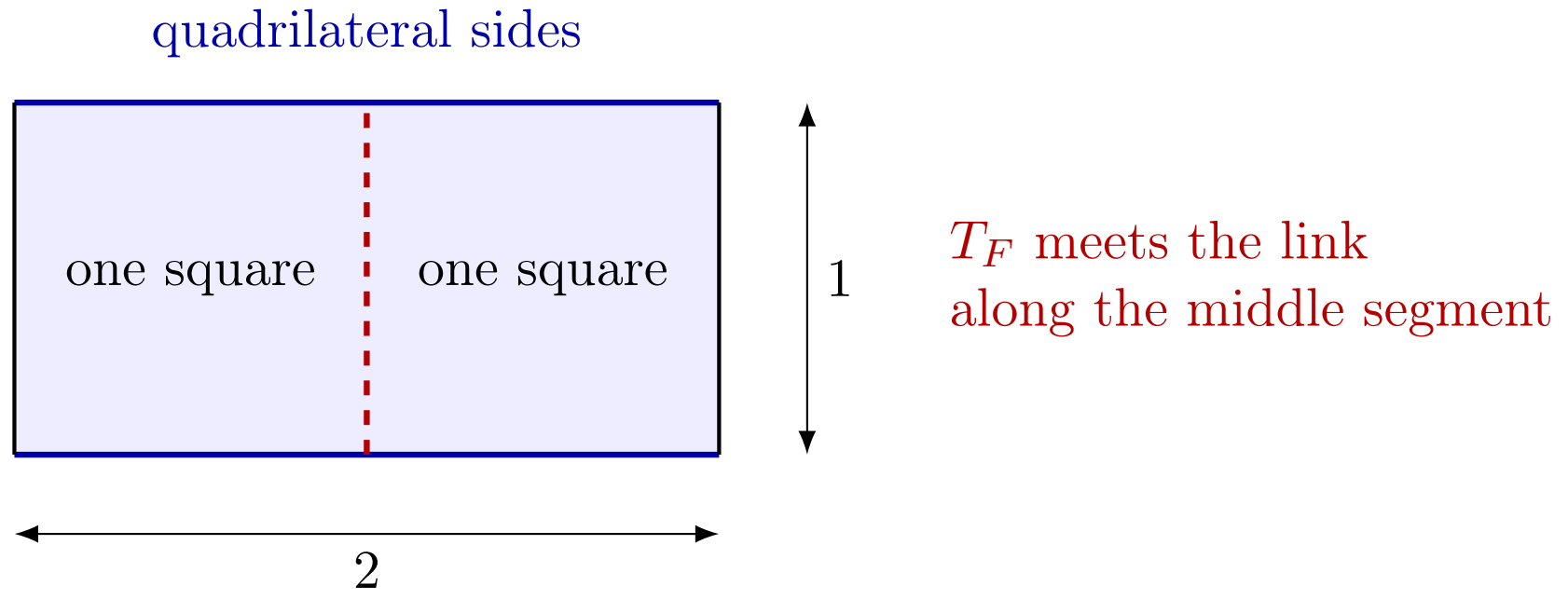

The convex hull of the three nonsquare cusps is an ideal triangle
$T_F$ whose relative interior lies in the relative interior of a
doubled facet $F$. This is its
original gluing triangle. Each side of $T_F$ is the interior diagonal
of a quadrilateral face of $F$. At a nonsquare cusp, $T_F$ meets
the horospherical rectangle along the transverse segment halfway
along its long direction. The rectangle and transverse segment
are shown in Figure~\ref{fig:rectangle}.

At a square cusp of $\DBL$, exactly one incident polygonal face
is a quadrilateral. At a cusp of $\OCT$, none is. These facts will be used both in the equality reduction and in the
general octahedral cutting argument.

\section{The facet types at the second volume}\label{sec:equalityspectrum}
In this section alone assume $\Vol(P)=8\pi^2/3$. Then $W=12$ and
$\sum_Fd(F)=24$. A facet is \emph{positive} when its defect is
positive; all others are regular ideal octahedra. Let $N$ be the
number of positive facets and $m$ their minimum defect. For any
positive facet, Lemma~\ref{lem:injection} gives
\begin{equation}\label{eq:equalitybudget}
 24\ge d(F)+m\bigl(q(F)+s(F)\bigr).
\end{equation}
We shall prove the following statement before making any use of
Theorem~\ref{thm:chain}.
\begin{proposition}\label{prop:equalityspectrum}
There are exactly eight positive facets, each isometric to $\DBL$.
All other facets are isometric to $\OCT$.
\end{proposition}
The three quadrilateral faces of $\DBL$ have no common edges and
meet pairwise at its three nonsquare cusps. Every triangle contains
at most one of those cusps; the quadrilaterals contain two each,
as opposite vertices. At a square cusp exactly one incident face
is quadrilateral. These are the intrinsic facts of
Section~\ref{sec:doubled}.

\subsection{Acute intersections}

We first explain why pairwise facet adjacency determines an actual
common face in the situation used below. Give each supporting wall
an outward unit spacelike normal. Distinct wall products are
nonpositive for an acute-angled Coxeter chamber
\cite[p.~42, \S2.1]{Vinberg}. Let $n_1,\dots,n_k$ be mutually
orthogonal actual wall normals, and let $p$ be a strict future
timelike interior point. Set
\[
 u=p-\sum_{i=1}^k\ip{p}{n_i}n_i.
\]
Orthogonality gives
\[
 u^2=p^2-\sum_{i=1}^k\ip{p}{n_i}^2<0,
 \qquad \ip{u}{n_i}=0.
\]
For every omitted actual wall normal $m$ we also have
\[
 \ip{u}{m}
 =\ip{p}{m}-\sum_{i=1}^k\ip{p}{n_i}\ip{n_i}{m}
 \le\ip{p}{m}<0.
\]
Indeed, both factors in each product being subtracted are
nonpositive. The vector $u$ is future timelike, since
$\ip{p}{u}=u^2<0$. After normalization it is therefore an actual
point of the selected codimension-$k$ face, strict against every
other wall.

In particular, three pairwise adjacent right-angled facets of $P$
meet along an actual edge. Four pairwise adjacent facets would meet
at a finite vertex, which is impossible because $P$ is entirely
ideal. If $G$ and $H$ are neighbors of $F$ and are adjacent to each
other anywhere in $P$, their two faces in $F$ consequently share
an edge. Thus the intrinsic facet picture also controls these
possible adjacencies in the ambient polytope.

\subsection{The minimum positive defect is three}\label{first-reduction-the-minimum-positive-defect-is-three}

We eliminate all possible values of $m$ other than three.

\paragraph{Antiprisms of defect at least four.}
Suppose that a facet is $A_k$ with $k\ge5$. Its defect is
$d=2k-6\ge4$, and $q=2$, $s=2k=d+6$. Applying the universal
bound $m\ge2$ in \eqref{eq:equalitybudget} gives
\[
 24\ge d+2(d+8)=3d+16\ge28,
\]
a contradiction. Thus no such antiprism occurs.

\paragraph{Defect two.}
A defect-two facet $F$ would be the square antiprism $A_4$.
It has eight nonsquare cusps and two quadrilateral faces; each cusp
belongs to exactly one of these quadrilaterals. If one face of a
rectangular parallelepiped is nonsquare, at least one of the other two
opposite-face pairs is nonsquare. Otherwise the corresponding two
side ratios would both be one, forcing all three sides to be equal.
At each cusp of $F$, at least two of its four ridge-neighbor facets
are consequently nonoctahedral. At most one is the quadrilateral
neighbor through that cusp, so there is a positive triangular
neighbor.

Covering all eight cusps requires at least three distinct triangular
neighbors, since each triangular face contains only three cusps.
They are distinct from $F$, its two quadrilateral neighbors, and its
eight nonsquare-cusp opposites. The last distinctness follows from
the same asymptotic-wall argument as in
Lemma~\ref{lem:injection}: a cusp opposite cannot be any ridge
neighbor. There are therefore at least
$1+2+8+3=14$ positive facets. Their total defect is at least $28$,
contrary to $\sum_Fd(F)=24$. Hence $m\ge3$.

\paragraph{Minimum defect at least five.}
If $m\ge5$, there are at most four positive facets. None is an
antiprism by the preceding exclusion, so each has at least three
nontriangular faces by Lemma~\ref{lem:two}. The neighbors across
those faces are distinct and positive. There must therefore be
exactly four positive facets, all pairwise adjacent. They would
meet at a finite vertex, a contradiction.

\paragraph{Minimum defect four.}
If $m=4$, there are at most six positive facets. Choose one,
say $F$, with defect four. A pentagonal face would force $A_5$ by
the antiprism equality case, already excluded. The published
ten-vertex classification therefore identifies $F$ with the
octahedrite in \cite[Figure~7]{Kolpakov}. Its four quadrilaterals
form two edge-adjacent pairs.

Let $G_1,\dots,G_4$ be the positive neighbors across these four
quadrilaterals. If two of them are adjacent, they belong to one of
the prescribed pairs, by the acute-intersection argument. Thus,
besides $F$ and at most one partner, each $G_i$ has only one other
possible positive neighbor.

Each $G_i$ needs at least three neighbors across nontriangular
faces. The additional facet must therefore exist, and $q(G_i)=3$.
A defect-four facet of the allowed type
has four nontriangular faces, so none of the $G_i$ has defect four.
Each has defect at least five. Including $F$ and the necessary
further facet gives a total defect at least
\[
 4+4\cdot5+4=28,
\]
again a contradiction. The only remaining possibility is $m=3$;
in particular some facet is $\DBL$.

\subsection{Exactly eight positive facets}\label{exactly-eight-positive-facets-all-d9}

A facet of type $\DBL$, its three quadrilateral neighbors, and its
three nonsquare-cusp opposites give seven positive facets. Since
$m=3$ and the total defect is $24$,
\[
 7\le N\le8.
\]
If $N=8$, all eight defects are three and the proposition follows.
Suppose, for a contradiction, that $N=7$. For every doubled facet,
the seven facets just counted exhaust all positive facets. Thus
\begin{equation}\label{eq:seven-triangular}
 \text{no triangular ridge of a doubled facet has a positive neighbor.}
\end{equation}
There are only three ways to distribute the excess defect above
seven copies of three: one defect six; one defect five and one
defect four; or three defects four. In particular at most three
facets have defect greater than three.

We first exclude defect four. Let $F$ have that defect, with its
four quadrilaterals paired by edge adjacency. If a quadrilateral
neighbor $G$ were doubled, let $H$ be the neighbor across the paired
quadrilateral of $F$. The shared edge of these two quadrilaterals
belongs to $F,G,H$. Inside $G$, the faces cut out by $F$ and $H$
share this edge. The first is a quadrilateral. The second cannot
be a quadrilateral, because two quadrilaterals of a doubled facet
never share an edge; it cannot be a triangle either, by
\eqref{eq:seven-triangular}, since $H$ is positive. This is
impossible. All four quadrilateral neighbors of $F$ would therefore
have defect greater than three, but there are at most two such
facets other than $F$.

The only remaining partition is one defect-six facet $F$ and six
doubled facets. Every nontriangle of $F$ is shared with a doubled
facet and is therefore a quadrilateral. There are six of them,
and no two share an edge: otherwise the same three-facet argument
would force either edge-adjacent quadrilaterals or a positive
triangular neighbor in a doubled facet.

The facet $F$ has twelve vertices, and its six quadrilaterals give
$24$ quadrilateral--vertex incidences. At a four-valent vertex at
most two incident faces can be quadrilaterals, since consecutive
quadrilaterals would share an edge. The incidence count forces
exactly two at every vertex, alternating with two triangles.
Make a graph $J$ with one vertex for each quadrilateral and one
edge for each original vertex of $F$, joining its two incident
quadrilaterals. It is a simple four-regular graph on six vertices.
Indeed, two distinct quadrilateral faces cannot share two vertices
without sharing their joining edge: their intersection is a face
of the convex three-polyhedron. The complement of $J$ is a perfect
matching, so $J$ is the octahedral graph.

We also need its embedding to recover $F$. Draw each new edge
from a center of one quadrilateral through the common original
vertex to the center of the other. Alternation at every original
vertex gives a spherical embedding whose eight triangular regions
are the eight original triangles. Their boundary cycles are
distinct. A connected embedded graph with more than three vertices
cannot have both sides of the same three-cycle as faces, since
there would be no side left for another vertex or edge. The
octahedral graph has exactly eight three-cycles, so all of them are
the face boundaries just found. Recovering $F$ puts one vertex on
each edge of $J$, one quadrilateral around each vertex of $J$, and
one triangle in each of its eight triangular regions. This is the
cuboctahedral face lattice. The argument determines its actual
face incidences as well as the abstract graph.

All twelve cusp rectangles of the ideal right-angled cuboctahedron
are nonsquare, with aspect ratio $\sqrt2$. For an exact model, take
Lorentz coordinates with ideal vertices
\[
 (1,x),\qquad
 x\text{ a permutation of }(\pm1,\pm1,0)/\sqrt2.
\]
Its six square-wall normals are $(1,\pm\sqrt2 e_i)$, of norm one,
and its eight triangle-wall normals are
$(1,(\pm1,\pm1,\pm1)/\sqrt2)$, of squared norm $1/2$.
At the vertex $(1,1/\sqrt2,1/\sqrt2,0)$, the opposite square
normals sum to $(2,\sqrt2,\sqrt2,0)$, and the opposite
\emph{unnormalized} triangle normals have the same sum. Passing to
unit triangle normals multiplies the latter sum by $\sqrt2$.
The horospherical width ratio is therefore $\sqrt2$, and symmetry
gives this ratio at every vertex. In scaled Klein coordinates the
model is the cube $|X_i|\le1$ truncated by
$\sum|X_i|\le2$, which also makes its face structure explicit.
Three-dimensional right-angled rigidity identifies $F$ with this
model.

Applying Lemma~\ref{lem:injection} to $F$ now gives
$1+6+12=19$ positive facets, contradicting $N=7$. Thus $N=8$,
and every positive facet has defect three and is isometric to
$\DBL$. All others are regular ideal octahedra. This proves
Proposition~\ref{prop:equalityspectrum}. There are twelve
quadrilateral ridges, and every other ridge is triangular.

\section{Cusp geometry under the octahedral facet hypothesis}\label{sec:chainproof}
We now assume \eqref{eq:spectrum}, allowing arbitrary volume, and
prove Theorem~\ref{thm:chain} from the intrinsic geometry of the facets.

\begin{lemma}\label{lem:box}
Every noncubic cusp section of $P$ is a rectangular parallelepiped
with side lengths in the ratio $2:1:1$.
Its four nonsquare faces belong to doubled facets and form a belt;
its two square end faces belong to ordinary octahedral facets.
Every adjacency within the belt is a quadrilateral ridge of $P$.
\end{lemma}
\begin{proof}
Each face of a cusp section is intrinsic to a facet of $P$.
Its aspect ratio is therefore one or two. If the three side lengths
were distinct, the ratio of the largest to the smallest would exceed
two, since the two consecutive unequal ratios must both be two.
Thus a noncubic parallelepiped has proportions $2:1:1$ or $2:2:1$.

Consider a $2:2:1$ parallelepiped. The four nonsquare faces must be doubled
facets. In each, the two long sides are quadrilateral sides. Each
square end face consequently has four incident quadrilateral sides.
An octahedral cusp has none, and a square cusp of a doubled facet
has exactly one. Neither permitted facet type can supply an end
face. The $2:2:1$ parallelepiped is impossible.

In a $2:1:1$ parallelepiped, the four long rectangles form the belt. Their
mutual sides are the long sides and hence are quadrilateral ridges.
The square ends meet the belt only along its short, triangular
sides. All four incidences at an end are triangular, so the end
cannot be a square cusp of a doubled facet. It is octahedral.
\end{proof}

\section{The internal triangles form closed surfaces}

Join two doubled facets when they share a quadrilateral ridge.
The distinguished triangles in these facets have the same internal
diagonal in that ridge, so we glue them along this diagonal. The
gluing identifies the common diagonal pointwise.

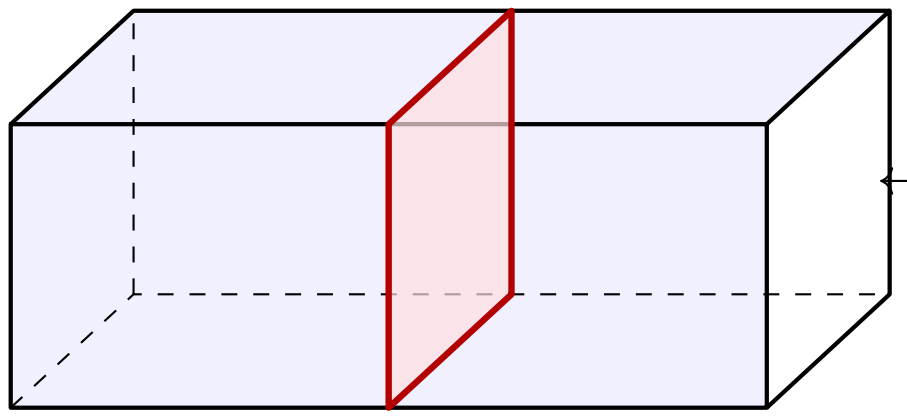
\begin{figure}[htbp]
\centering
\begin{tikzpicture}[scale=1.25,line join=round]
 \coordinate (A) at (0,0); \coordinate (B) at (4,0);
 \coordinate (C) at (4,1.5); \coordinate (D) at (0,1.5);
 \coordinate (E) at (.65,.6); \coordinate (F) at (4.65,.6);
 \coordinate (G) at (4.65,2.1); \coordinate (H) at (.65,2.1);
 \fill[blue!6] (A)--(B)--(C)--(G)--(H)--(D)--cycle;
 \draw[thick] (A)--(B)--(C)--(D)--cycle (B)--(F)--(G)--(C)
 (D)--(H)--(G);
 \draw[dashed] (A)--(E)--(F) (E)--(H);
 \fill[red!12,opacity=.7] (2,0)--(2.65,.6)--(2.65,2.1)--(2,1.5)--cycle;
 \draw[very thick,red!70!black] (2,0)--(2.65,.6)--(2.65,2.1)--(2,1.5)--cycle;
 \node[below] at (2,-.15) {middle of the long direction};
 \draw[->] (5.5,1.2)--(4.6,1.2);
 \node[anchor=west,align=left] at (5.55,1.2)
 {square end:\\ordinary octahedron};
\end{tikzpicture}
\caption{Four internal triangles meet a $2:1:1$ cusp section
in its transverse middle square. Its boundary is the four-segment
link used to prove that the compactified triangle complex is a
two-manifold. The filled square indicates the transverse section;
only its boundary lies on the surface of the parallelepiped.}
\label{fig:belt}
\end{figure}

\begin{lemma}\label{lem:surface}
After adjoining their ideal vertices, the glued triangles form a
finite closed two-manifold. Its connected components are exactly
the components of the quadrilateral adjacency graph of doubled facets.
Every edge lies in two triangles and every vertex lies in four.
\end{lemma}
\begin{proof}
The interior of a distinguished triangle lies in the relative
interior of its facet. Interiors of different facets are disjoint.
An edge of the triangle is in the relative interior of an actual
quadrilateral ridge, which lies in precisely two facets. Thus each
edge has two incident triangles and has a disk neighborhood at
each of its finite points.

It remains to check the ideal vertices. By Lemma~\ref{lem:box},
the triangles incident to one such point come from the four doubled
facets in a single $2:1:1$ belt. Their four middle segments join
end to end around the parallelepiped, giving the simple square shown in
Figure~\ref{fig:belt}. This square is the link of the point in the glued surface.
Coning the link gives a disk after the ideal point is adjoined.
There are no other triangle vertices.

All four doubled facets at such a vertex are already connected
through their quadrilateral belt. Therefore identifying ideal
vertices does not join two different components of the quadrilateral
adjacency graph. The components of the graph and of the surface
coincide.
\end{proof}

\begin{lemma}\label{lem:euler}
A connected component containing $n$ triangles contains exactly
eight triangles, twelve edges and six vertices.
\end{lemma}
\begin{proof}
Write $e$ and $v$ for its edge and vertex counts. The incidence
counts of Lemma~\ref{lem:surface} give
\[
 e=\frac{3n}{2},\qquad v=\frac{3n}{4},\qquad
 \chi=v-e+n=\frac n4.
\]
For a connected closed cellulated surface, $\chi\le2$.
We give a direct proof.
Choose a spanning tree $A$ of its one-skeleton. Delete from the
dual graph the edges dual to $A$. The remaining dual graph is
connected. Otherwise, the boundary between one nonempty proper
set of triangular faces and its complement would be a nonempty
subgraph of $A$. At each vertex its degree is even: a walk around
the circular link enters and leaves the chosen set equally often.
This gives a nonempty Eulerian subgraph of a tree, which is
impossible. We may therefore choose a dual spanning tree disjoint
from $A$. The two trees use $v-1$ and $n-1$ distinct edges, so
$e\ge v+n-2$ and $\chi\le2$. It follows that $n\le8$, and the
integer $\chi=n/4$ also shows that $n$ is divisible by four.

We next obtain a lower bound within the same component. Choose
one doubled facet $F$. Its three quadrilateral faces have three
distinct doubled neighbors. Its three nonsquare cusps also have
three opposite doubled facets. The latter are distinct from one
another and from every ridge neighbor: their walls are asymptotic
to the wall of $F$ at exactly one specified ideal point, whereas
a ridge neighbor intersects its wall along a plane. These opposite
facets belong to the same component, because opposite belt members
are connected by the four-cycle of quadrilateral adjacencies.
Together with $F$ we have found seven distinct facets in the
component. Thus $7\le n\le8$ and $4\mid n$, forcing $n=8$.
The values of $e$ and $v$ follow.
\end{proof}

\section{From six cusps to eight wall normals}

Fix one component. Make a graph $J$ whose six vertices are its
cusps and whose edges are the twelve internal diagonals. Every
vertex has degree four.

\begin{lemma}\label{lem:graph}
The graph $J$ is the octahedral graph. Its eight triangles are the
eight distinguished facet triangles. The doubled facets can
therefore be indexed by $s\in\{\pm1\}^3$ so that quadrilateral
adjacency changes exactly one sign.
\end{lemma}
\begin{proof}
First, $J$ has no parallel edges. If two different quadrilateral
ridges joined the same two cusps, their incident facets would
include at least three of the four belt members at either cusp.
Among three belt members there are two opposite ones. Their walls
are asymptotic at that cusp and cannot contain the second cusp.
This is a contradiction. Two distinct ridges cannot have the same
pair of incident facets, since the intersection of two convex
facets is one face. Thus $J$ is simple.

A simple four-regular graph on six vertices has complement a
perfect matching, so it is the octahedral graph. Each distinguished
triangle gives a three-cycle in $J$. Two facets cannot give the
same three-cycle: its geometric triangle would then have an
interior point in the relative interiors of two different facets.
The octahedral graph has exactly eight three-cycles, obtained by
choosing one vertex from each of its three opposite pairs. All
eight occur. Label these opposite pairs by $1,2,3$, with a sign
on each member. A triangle, and hence a doubled facet, is indexed
by its three choices of sign. Two triangles share an edge exactly
when one sign changes.
\end{proof}

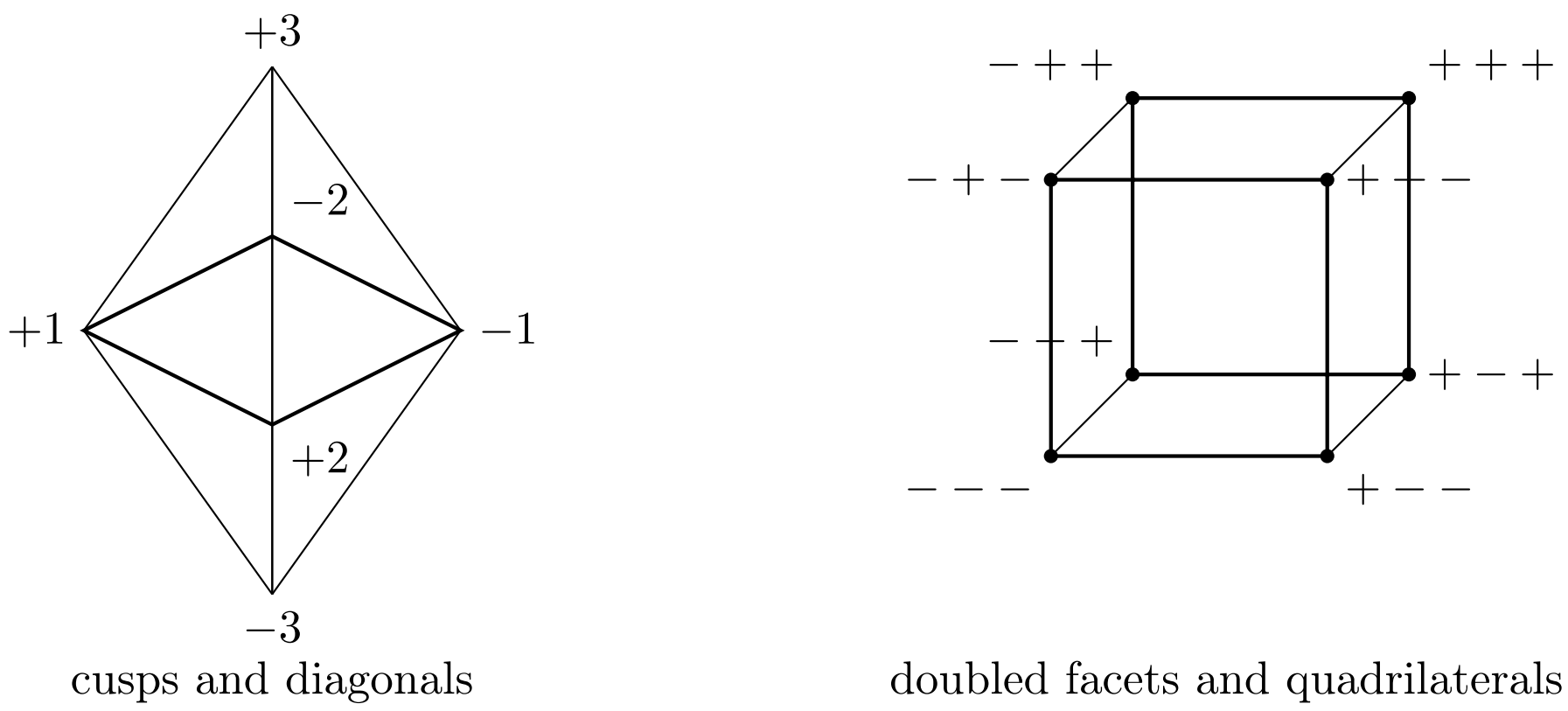
\begin{figure}[htbp]
\centering
\begin{tikzpicture}[scale=.95,every node/.style={font=\small}]
 \coordinate (xp) at (-1.5,0); \coordinate (yp) at (0,-.75);
 \coordinate (xm) at (1.5,0); \coordinate (ym) at (0,.75);
 \coordinate (zp) at (0,2.1); \coordinate (zm) at (0,-2.1);
 \draw[thick] (xp)--(yp)--(xm)--(ym)--cycle;
 \foreach \p in {xp,yp,xm,ym}{\draw (zp)--(\p) (zm)--(\p);}
 \node[left] at (xp) {$+1$}; \node[right] at (xm) {$-1$};
 \node[below right] at (yp) {$+2$}; \node[above right] at (ym) {$-2$};
 \node[above] at (zp) {$+3$}; \node[below] at (zm) {$-3$};
 \node at (0,-2.8) {cusps and diagonals};
 \begin{scope}[xshift=6.2cm,yshift=-1cm]
 \coordinate (a) at (0,0); \coordinate (b) at (2.2,0);
 \coordinate (c) at (2.2,2.2); \coordinate (d) at (0,2.2);
 \coordinate (e) at (.65,.65); \coordinate (f) at (2.85,.65);
 \coordinate (g) at (2.85,2.85); \coordinate (h) at (.65,2.85);
 \draw[thick] (a)--(b)--(c)--(d)--cycle (e)--(f)--(g)--(h)--cycle;
 \draw (a)--(e) (b)--(f) (c)--(g) (d)--(h);
 \foreach \p/\lab/\pos in {a/{---}/below left,b/{+--}/below right,
 c/{++-}/right,d/{-+-}/left,e/{--+}/above left,f/{+-+}/right,
 g/{+++}/above right,h/{-++}/above left}
 {\fill (\p) circle (1.6pt);\node[\pos] at (\p) {$\lab$};}
 \node at (1.4,-1.8) {doubled facets and quadrilaterals};
 \end{scope}
\end{tikzpicture}
\caption{The two incidence graphs. The six-cusp graph is octahedral;
its eight triangles index the doubled facets. The dual cube records
quadrilateral adjacency. The signs specify which cusp from each
opposite pair belongs to a distinguished triangle.}
\label{fig:graphs}
\end{figure}

\section{Equal widths put the triangles in one hyperplane}

We work in $\mathbb R^{4,1}$ with form
$\ip{x}{x}=-x_0^2+x_1^2+\cdots+x_4^2$. A wall has an outward
unit spacelike normal $n$, and its halfspace is $\ip{x}{n}\le0$.
Let $n_s$ be the normals of the eight doubled facets just indexed.

\begin{lemma}\label{lem:width}
On each square of the cube in Figure~\ref{fig:graphs}, the sums of
the two opposite pairs of normals agree.
\end{lemma}
\begin{proof}
The square is a cusp belt. Its opposite wall pairs bound the two
short directions of the $2:1:1$ parallelepiped, so their Euclidean separations
on a common horosphere agree. The use of unit normals matters.
With cusp null vector $\xi=(1,0,0,0,1)$, the vertical wall $x_i=a$
in upper-half-space has unit normal $(a,e_i,a)$. The opposite wall
$x_i=b$, with the other outward orientation, has normal
$(-b,-e_i,-b)$. Their sum is $(a-b)\xi$. Thus its coefficient is
the separation, with the same normalization for both pairs. Equal
widths give equal sums.
\end{proof}

\begin{lemma}\label{lem:normal}
There are mutually orthogonal vectors $a,b_1,b_2,b_3$ such that
\[
 n_s=a+s_1b_1+s_2b_2+s_3b_3,\qquad
 a^2=-\tfrac12,\quad b_i^2=\tfrac12.
\]
\end{lemma}
\begin{proof}
The equations of Lemma~\ref{lem:width} say that every mixed second
finite difference on the cube vanishes. Expanding $s\mapsto n_s$
in sign characters leaves only its constant and three linear
terms, giving the displayed affine expression. Expanding the
eight equations $n_s^2=1$ and comparing sign characters gives
\[
 \ip{a}{b_i}=0,\quad \ip{b_i}{b_j}=0\ (i\ne j),\quad
 a^2+b_1^2+b_2^2+b_3^2=1.
\]
Adjacent facets meet orthogonally. If $s'$ differs from $s$ only
in coordinate $i$, then
$0=\ip{n_s}{n_{s'}}=1-2b_i^2$. Thus every $b_i^2=1/2$, and
the remaining identity gives $a^2=-1/2$.
\end{proof}

The span $U$ of these four vectors is Lorentzian and has dimension
four. Its orthogonal complement is a spacelike line. The
intersection of $U$ with $\HH^4$ is a genuine hyperbolic
hyperplane $H$, orthogonal to all eight doubled-facet walls.
In orthonormal coordinates on $U$, the wall normals are
$(1,s_1,s_2,s_3)/\sqrt2$. The future orientation is determined by
the sums of opposite wall normals at a cusp, which are positive multiples of
the cusp vectors. Their halfspaces cut out
\begin{equation}\label{eq:section}
 S=\{(t,x)\in\HH^3:|x_1|+|x_2|+|x_3|\le t\},
\end{equation}
the regular ideal octahedron. Its six null rays are $a\pm b_i$.
They are exactly the actual cusps of the component, because the
opposite belt-normal sums are $2(a\pm b_i)$.

\begin{figure}[htbp]
\centering
\begin{tikzpicture}[scale=1.15]
 \draw[gray] (0,0) circle (1.75);
 \coordinate (a) at (0,1.75); \coordinate (b) at (1.75,0);
 \coordinate (c) at (0,-1.75); \coordinate (d) at (-1.75,0);
 \fill[blue!7] (a)--(b)--(c)--(d)--cycle;
 \draw[thick] (a)--(b)--(c)--(d)--cycle;
 \draw[very thick,red!70!black] (a)--(c);
 \fill[blue!65!black] (a) circle (2pt) (c) circle (2pt);
 \node[above] at (a) {nonsquare cusp};
 \node[below] at (c) {nonsquare cusp};
 \node[anchor=west,align=left] at (2.4,.3)
 {$S\cap R$ is the red diagonal.\\Its finite points lie in the\\
 relative interior of $R$.};
 \draw[->] (2.2,.1)--(.15,0);
\end{tikzpicture}
\caption{A quadrilateral ridge $R$ in the Klein disk. Its two
nonsquare cusps are opposite vertices. Their diagonal is an edge
of the cutting octahedron; its finite points lie in the relative
interior of $R$.}
\label{fig:diagonal}
\end{figure}
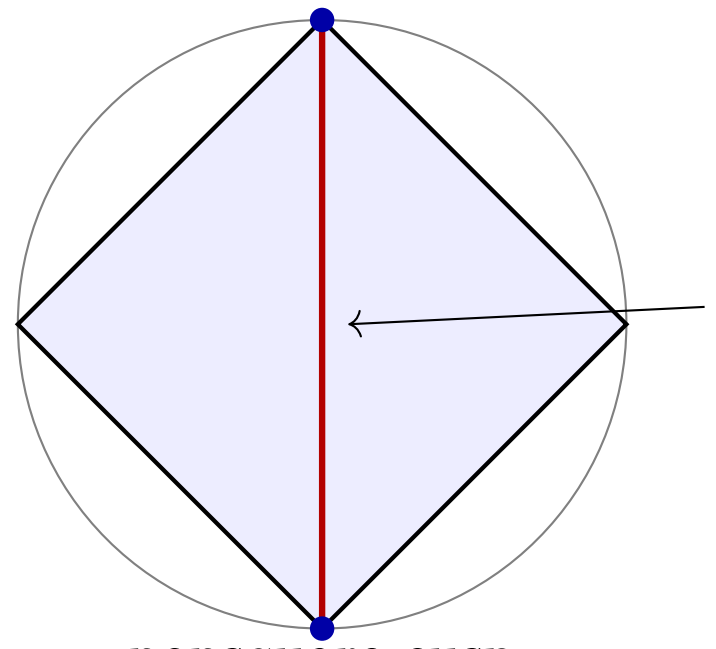

Since these six ideal points belong to the closure of $P$, their
hyperbolic convex hull $S$ is contained in $P$. Conversely the
eight actual facet inequalities imply $P\cap H\subset S$.
Therefore
\[
 P\cap H=S.
\]
Each face of $S$ is its facet's distinguished internal triangle.
Each edge of $S$ is the internal diagonal of an actual
quadrilateral ridge, as illustrated in Figure~\ref{fig:diagonal}.

\section{Cutting along the octahedron}

\begin{lemma}\label{lem:cut}
No facet outside the chosen eight meets $S$ at a finite point.
Cutting along $H$ gives two finite-volume entirely ideal
right-angled four-polytopes. The eight doubled facets become
ordinary octahedra; every other facet is uncut.
\end{lemma}
\begin{proof}
Let $F$ be any outside facet. Its supporting functional is
nonpositive on $S$. If it vanishes at a relative interior point
of a face of $S$, it vanishes on that face. A finite point of $S$
lies in its interior, in a triangular face interior, or in an
edge interior. A triangular face interior belongs to the relative
interior of one of the chosen doubled facets, so it cannot belong
to $F$. An edge interior lies in the relative interior of a ridge
whose only two incident facets are the chosen doubled facets.
It cannot belong to $F$ either. A supporting functional vanishing
at an interior point would contain all of $S$, contradicting the
face-interior case. Thus $F$ has no finite intersection with $S$.

If $F$ crossed $H$, a finite segment in $F$ would meet
$P\cap H=S$ at a finite point. Consequently no outside facet
crosses. The chosen eight facets are cut along their genuine
gluing triangles into octahedra. The center of $S$ is strict on
all eight wall inequalities and, by what was just proved, on
every other wall inequality. It is an interior point of $P$;
both halves have nonempty interior.

In the Klein model, a vertex of a half lying on $H$ must be a
vertex of $S$. Otherwise it lies in a nontrivial segment of $S$
and is not extreme. Away from $H$, a vertex is an old vertex of
$P$. The only vertices of $S$ are its six old ideal cusps, so
no finite vertex is introduced. The new facet $S$ meets only the
eight cut walls, all orthogonal to $H$. Every new dihedral angle
is right; all other angles are inherited. Finite volume follows
by containment in $P$.
\end{proof}

The facet spectrum \eqref{eq:spectrum} is preserved by the cut.
There are exactly eight fewer doubled facets, counted over both
pieces. Repeating the construction therefore terminates.

There is also a useful geometric disjointness statement. Every
outside doubled facet lies on one side of the first cut. Two
quadrilateral-adjacent outside facets cannot lie on opposite
sides, since their shared ridge would give a finite intersection
with $S$. Thus a whole other component lies on one side. Its
section, being the convex hull of its six cusps, lies on that
side as well. A noncubic cusp belongs to only one belt component
and has octahedral ends; two sections cannot share an ideal
vertex. The outside-facet argument excludes finite boundary
contact. Hence the sections are pairwise disjoint, including
their ideal endpoints.

\section{Recognizing a piece with only octahedral facets}

\begin{lemma}\label{lem:recognition}
An entirely ideal right-angled four-polytope whose facets are all
isometric to $\OCT$ is isometric to $\CELL$.
\end{lemma}
\begin{proof}
Place one facet $F$ on a facet of a model regular 24-cell. Its
eight triangular faces agree with those of the model facet.
For any one of these triangles, the adjacent supporting wall is
uniquely determined by the triangle and orthogonality to $F$'s
wall: its outward normal is the intrinsic outward unit normal
of that triangle within $F$'s wall.

In this adjacent wall, its octahedral facet is determined by the
shared ideal triangle and the side on which it lies. To justify
this last assertion, take any isometry between two regular ideal
octahedra with that triangular face. Every isometry of the face
extends to a symmetry of the octahedron, so the isometry can be
adjusted to fix the face vertices pointwise. An isometry of
$\HH^3$ fixing three distinct ideal points is the identity or
reflection in their plane; the specified side selects the
identity. The neighboring facet therefore agrees with the model
octahedron.

Repeat across its faces. The facet adjacency graph of the model
24-cell is connected. The propagation matches all 24 model
facets, with their full vertex sets, to actual facets of the
piece. Their halfspaces put the piece inside the model cell.
Conversely the closure of the piece contains every ideal model
vertex, and convexity contains their hull, the model cell.
The two polytopes agree.
\end{proof}

Each final piece is consequently a 24-cell. If the original
polytope has $N$ doubled facets, Lemma~\ref{lem:euler} gives
$N=8k$. Each of the $k$ cuts increases the number of pieces by
one, so there are $k+1$ cells. Their adjacency graph is a tree:
starting with one piece, each separating cut replaces one vertex
of the current adjacency tree by two vertices joined by one edge.

\section{Why the tree is a chain}

The restriction on facet types has one more consequence. In the
Lorentz model a regular 24-cell has outward unit normals
\[
 (1,r),\qquad
 r\in\{(\pm1,\pm1,0,0)\text{ and all coordinate permutations}\}.
\]
The spatial vectors are the roots of $D_4$. Distinct roots have
dot products $-2,-1,0,1$, with $-2$ only for opposite roots.

\begin{lemma}\label{lem:opposite}
If a cell in the decomposition meets two cutting sections,
these are opposite facets of that cell.
\end{lemma}
\begin{proof}
The sections have disjoint closures. If their roots are $r_1,r_2$,
then $r_1\cdot r_2$ is negative: dot product one gives an actual
common triangular ridge, and dot product zero gives a common
ideal vertex. Thus the dot product is $-1$ or $-2$.

Suppose it is $-1$. The vector $r_1+r_2$ is another $D_4$ root,
and its cell facet $Q$ is orthogonal to both cut walls, since
\[
 -1+(r_1+r_2)\cdot r_i=0\qquad(i=1,2).
\]
It cannot itself be a cutting facet, since it meets the two
sections along triangles. It is part of an original facet of $P$.
Across each of the two cuts there is a neighboring regular
24-cell. Its corresponding facet is the reflection of $Q$ and
lies in the same original supporting wall, because that wall is
orthogonal to the cut. This neighboring facet cannot be another
internal cut: it meets the first cut along a triangle, whereas
all cuts have disjoint closures. Thus the original facet contains
three regular ideal octahedra.

Their relative interiors are disjoint. If two had an open
overlap, a point of that overlap avoiding every other cell wall
would have an inward normal collar contained in both cells.
Their four-dimensional interiors would overlap, which is
impossible for pieces obtained by cuts. The original facet has
intrinsic volume at least three times $\Vol(\OCT)$; it can be
neither $\OCT$ nor $\DBL$. This contradicts
\eqref{eq:spectrum}. Hence the dot product is $-2$, and the two
facets are opposite.
\end{proof}

No cell can have three pairwise opposite facets. The adjacency
tree has maximum degree two and is a path. Once an end cell and
its first gluing facet are placed, every next cell is its
reflection across the gluing facet, and the next gluing facet in
an intermediate cell is uniquely the opposite one.

Every isometry of an octahedral facet extends to a symmetry of
the regular 24-cell. Indeed, the cell is flag-transitive. Given an isometry of an
octahedral facet, send one facet flag to its image by a cell
symmetry. This symmetry stabilizes the facet and has the required
restriction, because an octahedral isometry fixing a flag is the
identity. Thus every facet symmetry extends to the cell.
Two cells on opposite sides of the same facet are therefore
uniquely placed up to the ambient reflection. This proves the
isometry assertion in Theorem~\ref{thm:chain}.

\begin{figure}[htbp]
\centering
\begin{tikzpicture}[>=Latex,scale=1.1]
 \foreach \i in {0,1,2,3}{
  \draw[rounded corners,fill=blue!5] (2.2*\i,0) rectangle +(1.7,1.2);
  \node at (2.2*\i+.85,.6) {$\CELL_{\the\numexpr\i+1\relax}$};}
 \foreach \i in {0,1,2}{\draw[thick] (2.2*\i+1.7,.6)--(2.2*\i+2.2,.6);}
 \node[align=center] at (4.1,-.7)
 {two cuts in an intermediate cell\\must use opposite octahedral facets};
\end{tikzpicture}
\caption{The adjacency graph of the $24$-cells is a path.
The two cutting facets in an intermediate cell are opposite
by Lemma~\ref{lem:opposite}.}
\label{fig:chain}
\end{figure}
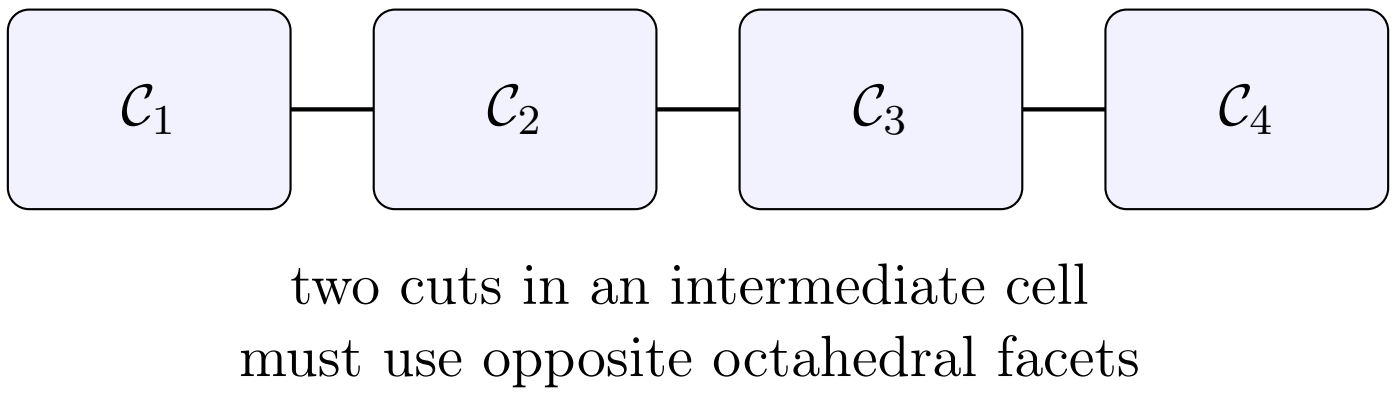

\section{Volumes and face numbers}

The regular ideal 24-cell has volume $4\pi^2/3$.
Volume additivity proves the formula in Theorem~\ref{thm:chain}. The common defect identities also give
$2W=3N=24k$, consistently with \eqref{eq:discrete}.

If there are $m=k+1$ cells, the complete face counts are
\[
 (v,e,r,f)=(18m+6,\ 72m+24,\ 68m+28,\ 14m+10).
\]
Indeed each gluing identifies six ideal vertices and removes two
facets while merging eight neighboring pairs. There are no extra
vertex identifications between nonadjacent cells: a common cusp
would lie on every separating cut along their path, whereas
different sections have no common ideal vertex. The edge and
ridge formulas then follow from the incidence identities above.
There are $8(m-1)$ doubled facets, $6m+18$ ordinary facets and
$12(m-1)$ quadrilateral ridges. All other ridges are triangles.
For $m=2$ this recovers $(42,168,164,38)$.

The facet hypothesis~\eqref{eq:spectrum} is used both to obtain the
cusp proportions of Lemma~\ref{lem:box} and to define the internal
triangles. The decomposition and volume formula apply to this class.

\section{Equality and attainment}\label{sec:attainment}
\begin{proof}[Equality in Theorem~\ref{thm:gap}]
If $\Vol(P)=8\pi^2/3$, Proposition~\ref{prop:equalityspectrum}
gives \eqref{eq:spectrum} and $N=8$. Theorem~\ref{thm:chain}
therefore gives $k=1$: there are exactly two regular 24-cells joined
across one octahedral facet. The uniqueness argument in
the chain proof shows that they are mutual reflections in that
facet. Conversely their reflection double is an entirely ideal
right-angled polytope of twice the 24-cell volume, as the explicit
construction below verifies. This proves the equality assertion.
\end{proof}

\subsection{An exact model of the attaining double}

In Lorentz form $(-1,1,1,1,1)$, take the 24 unit normals $(1,a)$,
where $a$ runs through the $D_4$ roots having exactly two nonzero
entries, each $\pm1$. In the time-one Klein section their ideal
vertices are the eight $\pm e_i$ and the sixteen
$(\pm\tfrac12,\pm\tfrac12,\pm\tfrac12,\pm\tfrac12)$.
The inequalities are $|x_i|+|x_j|\le1$; they are bounded and their
origin is strict. Here is a finite exact verification of the vertex list. For each of
the $\binom{24}{4}=10,626$ four-element subsets of the walls, form
the $4\times4$ matrix of their spatial coefficients. If its rank is
four, solve the four affine equalities over $\mathbb Q$ and retain
the solution exactly when it satisfies all $24$ inequalities.
Discard repeated solutions. Singular subsets contribute no candidate.
This procedure terminates because its loops and all its rational
linear systems have fixed finite size.

It is complete for this intersection. At a vertex of a full-dimensional
polytope, the active spatial wall coefficients span $\mathbb R^4$:
otherwise a nonzero vector orthogonal to them permits sufficiently
small perturbations in both directions while every inactive inequality
remains strict, contradicting extremality. Some four active walls
therefore have rank four and occur in the enumeration. Conversely a
feasible solution of four independent supporting equalities is a
vertex, since any segment through it in the intersection would have
both endpoints in those four equalities. The intersection here is
bounded, as follows directly from $|x_i|+|x_j|\le1$, and contains
the origin strictly. The exact calculation returns precisely the
eight signed coordinate vertices and the sixteen half-sign vertices
displayed above. Each has Euclidean norm one, so all are ideal.

Reflect across $n=(1,1,1,0,0)$ by
$x\mapsto x-2\ip{x}{n}n$. Remove the gluing wall from both cells
and identify coincident walls. The resulting 38 halfspaces contain
the union. Indeed, for an old wall normal $m\ne n$ and $x\in\CELL$,
\[
 \ip{x}{s_n(m)}=\ip{x}{m}-2\ip{m}{n}\ip{x}{n}\le0,
\]
since both $\ip{m}{n}$ and $\ip{x}{n}$ are nonpositive. The
mirrored containment follows by reflection. Conversely, on each side of the gluing hyperplane, their
inequalities together with that side's gluing inequality include all
inequalities of the corresponding original cell. Thus the new
intersection equals the union.

The interiors are disjoint, so the volume is twice $4\pi^2/3$.
In the Klein model this convex union is the convex hull
of the two finite vertex sets. Reflection and normalization give
$42$ distinct points on the ideal sphere, with six shared points.
Strict convexity of the sphere makes every one of these points an
extreme point of their hull. Hence this is the complete vertex set.

To reconstruct faces, first record, for each wall, the subset of these
vertices on its equality hyperplane. Start with the full vertex set
and repeatedly intersect known sets with these facet sets, discarding
the empty set and repetitions, until no new set appears. There are
at most $2^{42}$ subsets, so the procedure halts. Every generated set
is the vertex set of a face, and every face of a polytope is an
intersection of its containing facets. The closure is therefore
complete. Its dimension is the affine rank of its vertices, computed
over $\mathbb Q$. The resulting exact incidence counts are those
in the following table.
\begin{center}
\begin{tabular}{lrr}
\toprule Quantity&24-cell&Double\\\midrule
Vertices&24&42\\ Edges&96&168\\ Ridges&96&164\\ Facets&24&38\\
Triangular ridges&96&152\\ Quadrilateral ridges&0&12\\
Octahedral facets&24&30\\ Doubled-octahedron facets&0&8\\
Ridge defect $W$&0&12\\ Reflection-orbifold Euler characteristic&1&2\\
\bottomrule
\end{tabular}
\end{center}
Each of the eight doubled-octahedron facets has nine vertices, eighteen
edges, eight triangles and three quadrilaterals. There are six shared
cusps and 36 other cusps. Every cusp has six incident facets, eight
edges and twelve ridges. The facet--vertex check is
$30\cdot6+8\cdot9=6\cdot42$. This completes attainment and
the construction of the equality model.

\appendix
\section{Finite computations}\label{app:reproduction}

\subsection{Models and algorithms}

The supplementary computations verify the incidence data for the
$24$-cell and its double, the small equality-case models, and the
$D_4$ root relations used in Lemma~\ref{lem:opposite}.

The finite checks use only explicitly bounded sets. The original
$24$-cell test is the rational active-wall algorithm just proved;
the face test is the terminating facet-intersection closure. For the
cuboctahedral model, rescale the spatial coordinates so the Lorentz
form is $-x_0^2+(x_1^2+x_2^2+x_3^2)/2$. Its normals are then
$(1,\pm2e_i)$ and $(1,\pm1,\pm1,\pm1)$. The same vertex argument
runs over the $\binom{14}{3}=364$ triples of these fourteen rational
wall equations. For the six-cusp graph check, a
four-regular graph on six vertices has a perfect matching as its
complement. Enumerate the $\binom{15}{3}$ triples of possible missing
edges and retain those covering each vertex once. These are exactly
the fifteen labelled matchings; each complement has the eight
triangles and cube dual described in the proof. The cube-normal
check is rational row reduction of its six square equations in
eight unknowns, with kernel spanned by the constant and three sign
characters. These procedures halt by their explicit finite bounds.
Their completeness concerns these specified models and graph classes.

The checks give the two face vectors
$(24,96,96,24)$ and $(42,168,164,38)$, the cuboctahedron's twelve
cusps with squared aspect ratio two, the fifteen labelled six-cusp
graphs and cube-equation rank four. The separate $D_4$ check examines
all $\binom{24}{2}=276$ unordered root pairs: their spatial products
$-2,-1,0,1$ occur respectively $12,96,72,96$ times. It also checks
the $96$ common-neighbor cases and their $1,152$ reflected
octahedral-vertex incidences.

\subsection{Supplementary verification}

The Python code, finite-model certificates and verification records
are available in the accompanying repository~\cite{Repository}.
They implement the finite checks described above, including an
independent check of the saved model data.

Selected Lean~4 modules verify the internal logical deductions
for the volume gap, the equality reduction, the cut-count recursion
and the gluing implication. The cited classifications, finite-model
outcomes and geometric arguments constructing and recognizing the
cuts are supplied as explicit premises.

\section*{Acknowledgments}
The author used OpenAI's GPT-6-astra and GPT-6.1-sol as research and coding
assistants to investigate proof ideas, help strengthen the results, improve
the exposition, and generate and audit the Python and Lean~4 code
accompanying this paper. The author reviewed the resulting mathematics
and code and takes responsibility for the final content.

\end{document}